\documentclass{article}

\usepackage {amsmath, amssymb, amsbsy, amsthm, graphicx, verbatim, eufrak, eucal}

\newtheorem{theorem}{Theorem}
\newtheorem{lemma}{Lemma}

\newtheorem{definition}{Definition}
\newtheorem{corollary}{Corollary}
\newtheorem{example}{Example}
\newtheorem{remark}{Remark}

\numberwithin{theorem}{section}
\numberwithin{definition}{section}
\numberwithin{lemma}{section}
\numberwithin{corollary}{section}
\numberwithin{equation}{section}
\numberwithin{proposition}{section}
\numberwithin{example}{section}
\numberwithin{remark}{section}
\numberwithin{figure}{section}

\def\es{\varnothing}
\def\SB{\subseteq}

\def\aa{\alpha}
\def\bb{\beta}
\def\gg{\gamma}
\def\dd{\delta}

\def\bm{\boldsymbol m}
\def\bn{\boldsymbol n}
\def\bp{\boldsymbol p}
\def\bq{\boldsymbol q}

\def\bw{\boldsymbol w}

\def\b0{\boldsymbol 0}

\def\imp{\Rightarrow}
\def\EQ{\Longleftrightarrow}
\def\eq{\Leftrightarrow}

\def\FFF{{\cal F}}
\def\EEE{{\cal E}}
\def\CCC{{\cal C}}

\def\GGG{{\cal G}}
\def\HHH{{\cal H}}

\def\BBB{{\cal B}}

\def\WWW{{\cal W}}

\def\MMM{{\cal M}}

\def\PPP{{\cal P}}
\def\QQQ{{\cal Q}}

\def\SSS{{\cal S}}
\def\TTT{{\cal T}}
\def\UUU{{\cal U}}

\def\bmfL{\bold\,\mathfrak L\,}

\def\begeq{\begin{equation}}
\def\edeq{\end{equation}}

\def\roster{\begin{enumerate}}
\def\endroster{\end{enumerate}}

\begin{document}

\title{Fundamentals of Media Theory}
	\author{Sergei~Ovchinnikov\\ Mathematics Department\\San Francisco State University\\San Francisco, CA 94132\\sergei@sfsu.edu}

\date\today
\maketitle

\begin{abstract}\noindent
Media theory is a new branch of discrete applied mathematics originally developed in mid-nineties to deal with stochastic evolution of preference relations in political science and mathematical psychology. The theory focuses on a particular semigroup of `messages' acting as transformations of a set of `states', called a `medium', whose axioms are both strong and natural. The term `medium' stems from a particular application in which the transformations formalize the effects, on an individual, of `tokens' of information delivered by the environment---that is, the `medium'. However, many different types of examples can be found, ranging from learning spaces to hypercube computers, suggesting that this concept is ubiquitous. The paper presents very basic concepts and results of media theory and is aimed at a wide body of researchers in discrete applied mathematics.
\end{abstract}


\section{Introduction} \label{S:intro}

The term `media theory' was coined by Jean-Claude Falmagne in his founding paper~\cite{jF97} where basic concepts and results were introduced. That paper was followed by papers~\cite{sO00} and~\cite{jF02} in which the theory was further advanced.

A medium is an algebraic structure describing a mathematical, physical, or behavioral system as it evolves from one `state' to another, in a set of such states. Each set is characterized by a collection of binary features, and differs from some other neighbor state by just one of those features. This structure is formalized as an ordered pair $(\SSS,\TTT)$ consisting of a set $\SSS$ of states and a set $\TTT$ of tokens and specified by constraining axioms (see Section~\ref{S:axioms}). Tokens are transformations of the set $\SSS$; strings of tokens are messages of the medium. States, tokens, and messages are three fundamental notions of media theory. 

The set $\PPP$ of all partial orders on a given finite set $X$ is an example of a set of states that can be casted as medium. For any two distinct partial orders $P$ and $Q$, one can `walk' in $\PPP$ from $P$ to $Q$ by adding or removing a single ordered pair of elements of $X$. The transformations of $\PPP$ consisting in the addition or removal of some pair are tokens of the medium on $\PPP$. In terms of media theory, there is a concise message producing $Q$ from $S$. There are many other families of partial orders that can be casted as media, including linear orders, weak orders, semiorders, and interval orders~\cite{jF97,jF02,sO05}.  Additional examples of media include learning spaces~\cite{jD99} and hyperplane arrangements~\cite{sO05,sO06}.

Various stochastic applications of media theory have been made in the context of opinion polls and related situations~\cite{jF97,jF96,jF97a,jF97b,mR99}. Effective algorithms for visualization~\cite{dE05} and enumeration~\cite{dE06} of media have been developed.

The paper presents a concise introduction to basic concepts and results of media theory. Our exposition differs, in some details, from those given in~\cite{jF97,jF02} and~\cite{sO00}. Most notably, we use a system of just two constraining axioms instead of four original ones. This new system is equivalent to the old one but makes the underlying concepts more consistent. Naturally, all the results of media theory remain valid. This approach is also employed in the forthcoming monograph~\cite{dE07}. We also include some new results from graph theory (Section~\ref{S:bipartite}) that are crucial in establishing fine properties of messages and media.

In Sections~\ref{S:token systems} and~\ref{S:axioms} basic definitions and axioms are introduced and independence of the two constraining axioms is established. Some fundamental properties of tokens and messages are presented in Section~\ref{S:tokens}, which is preceeded by Section~\ref{S:example} where an important and in some sense `generic' example of a medium is given. Graphs are important tools in studying and representing media. Necessary facts from graph theory are presented in Section~\ref{S:bipartite}; graphs of media are introduced in Section~\ref{S:graphs}. We then proceed by establishing many properties of messages and media in Sections~\ref{S:circuits}--\ref{S:embeddings}. Graphs representing media are characterized as partial cubes and mediatic graphs in the last section of the paper.

\section{Token systems} \label{S:token systems}

Let $\SSS$ be a set of {\em states}. A {\em token} is a transformation $\tau:S\mapsto S\tau$. By definition, the identity function $\tau_0$  on $\SSS$ is not a token. Let $\TTT $ be a set of tokens. The pair $(\SSS,\TTT )$ is called a {\em token system}. To avoid trivialities, we assume that  $|\SSS|\geq 2$ and $\TTT\neq \es$.

Let $V$ and $S$ be two states of a token system $(\SSS,\TTT)$. Then $V$ is {\em adjacent} to $S$ if $S\neq V$ and $S\tau = V$ for some token $\tau\in\TTT$. A token $\tilde\tau\in\TTT$ is a {\em
reverse} of a token $\tau$ if for all distinct
$S,V\in\SSS$, we have 
$$
S\tau=V \quad\EQ\quad  V\tilde \tau=S.
$$
Two distinct states $S$ and $V$ are {\em adjacent} if $S$ is adjacent to $V$ and $V$ is adjacent to $S$.

\begin{remark}
{\rm In both examples of Figure~\ref{so-4-counter}, the state $V$ is adjacent to the state $S$, but these two states are not adjacent in either example.
}
\end{remark}

\begin{remark}
{\rm It is easy to verify that if a reverse of a token exists, then it is unique and the reverse of a reverse is the token itself; that is, $\tilde{\tilde\tau}=\tau$, provided that $\tilde\tau$ exists. In general, a token of a token system $(\SSS,\TTT)$ does not necessarily have a reverse in $(\SSS,\TTT)$. For instance, the token $\tau$ in Example [$\MMM$1] of Figure~\ref{so-4-counter} does not have a reverse in $\TTT$. It is also possible for a token to be the reverse of itself. For example, let $\SSS=\{S,V\}$, $\TTT=\{\tau\}$ where $\tau$ is the function defined by $S\tau=V$ and $V\tau=S$. Clearly, $\tilde\tau=\tau$.
}
\end{remark}

A {\em message} of a token system $(\SSS,\TTT)$ is a string of elements of the set $\TTT$. We write these strings in the form $\bm=\tau_1\tau_2\ldots\tau_n$. If a token $\tau$ occurs in the string $\tau_1\tau_2\ldots\tau_n$, we say that the message $\bm=\tau_1\tau_2\ldots\tau_n$ {\em contains} $\tau$. 

A message $\bm=\tau_1\tau_2\ldots\tau_n$ defines a transformation
$$
S\mapsto S\bm=((\ldots((S\tau_1)\tau_2)\ldots)\tau_n)
$$
of the set of states $\SSS$. By definition, the empty message defines the identity transformation $\tau_0$ of $\SSS$. If $V=S\bm$ for some message $\bm$ and states $S,V\in\SSS$, then we say that $\bm$ {\em produces} $V$ from $S$ or, equivalently, that $\bm$ {\em transforms} $S$ into $V$. More generally, if $\bm=\tau_1\ldots\tau_n$, then we say that $\bm$ {\em produces} a sequence of states $(S_i)$, where $S_0=S$ and $S_i=S\tau_1\ldots\tau_i$ for $1\le i\le n$.

If $\bm$ and $\bn$ are two messages, then $\bm\bn$ stands for the concatenation of the strings $\bm$ and $\bn$. We denote by $\widetilde\bm=\tilde\tau_n\ldots\tilde\tau_1$ the {\em reverse} of the message $\bm=\tau_1\ldots\tau_n$, provided that the tokens in $\widetilde\bm$ exist. If $\bn=\bm\bp\bm'$ is a message, with $\bm$ and $\bm'$ possibly empty messages, and $\bp$ non empty, then we say that $\bp$ is a {\em segment} of $\bn$.

The {\em content} of a message $\bm=\tau_1\ldots\tau_n$ is the set $\CCC(\bm)=\{\tau_1,\ldots,\tau_n\}$ of its distinct tokens. The content of the empty message is the empty set. We write $\ell(\bm)=n$ to denote the {\em length} of the message $\bm$ and assume that the length of the empty message is zero. It is clear that $|\CCC(\bm)|\le\ell(\bm)$ for any message $\bm$.

A message is {\em consistent} if it does not contain both a token and its reverse, and {\em inconsistent} otherwise. A message $\bm=\tau_1\ldots\tau_n$ is {\em vacuous} if the set of indices $\{1,\ldots,n\}$ can be partitioned into pairs $\{i,j\}$, such that $\tau_i$ and $\tau_j$ are mutual reverses.

A~message $\bm$ is {\em effective} (resp.~{\em ineffective})  for a state $S$ if $S\bm\neq S$ (resp. $S\bm=S$) for the corresponding transformation $\bm$. A~message $\bm=\tau_1\ldots\tau_n$ is {\em stepwise effective} for $S$ if  $S_k\neq S_{k-1}$, $1\leq k\leq n$, in the sequence of states produced by $\bm$ from $S$. A message is said to be {\em concise} for a state $S$ if it is stepwise effective for $S$, consistent, and any token occurs at most once in the message. 
A message is {\em closed} for a state $S$ if it is stepwise effective and ineffective for $S$.
When it is clear from the context which state is under consideration, we may drop a reference to that state.

If $\bm$ and $\bn$ are two concise messages such that $S\bm=V$ and $V\bn=S$ for some states $S$ and $V$, we call $\bm\bn$ a {\em $2$-gon} for $S$.

Some properties of the concepts introduced in this section are listed below. These properties are straightforward and will be used implicitly in this paper.

\roster
\item One must distinguish messages from transformations defined by these messages. For instance, for any token $\tau_i$, the two distinct messages $\bm=\tau_i\tau_i$ and $\bn=\tau_i$ of the token system displayed in Figure~\ref{example medium} define the same transformation of the set of states $\SSS$.

\item A consistent message may not contain a token which is identical to its reverse. Clearly, this also holds for concise messages.

\item The length of a vacuous message is an even number.

\item The reverse $\widetilde\bm$ of a concise message $\bm$ producing a state $V$ from a state $S$ is a concise message for $V$, provided that $\widetilde\bm$ exists.

\item Let $\bm=\tau_1\ldots\tau_n$ be a stepwise effective message for a state $S$. For any $i$, the state $S_{i+1}$ is adjacent to the state $S_i$ in the sequence of states produced by $\bm$. In general, there could be identical states in this sequence; a $2$-gon $\bm\widetilde\bm$ is an example (we assume that $\widetilde\bm$ exists).

\item Any segment of a concise message is a concise message for some state.

\item If $\bm$ is a concise message for some state, then $\ell(\bm)=|\CCC(\bm)|$.

\item A $2$-gon for a state $S$ is closed for $S$.
\endroster

\section{Axioms for a Medium} \label{S:axioms}

\begin{definition}\label{so-axioms} 
A token system $(\SSS,\TTT)$ is called a {\em medium} (on $\SSS$) if the following axioms are satisfied.
\begin{roster}
\item[]
\begin{roster}
\item[{[$\MMM$1]}] For any two distinct states $S$ and $V$ in $\SSS$ there is a
concise message transforming $S$ into $V$.
\item[{[$\MMM$2]}] A message which is closed for some state is vacuous.
\end{roster}
\end{roster}
A medium $(\SSS,\TTT)$ is {\em finite} if $\SSS$ is a finite set.
\end{definition}

{\begin{figure}[h!]
\centerline{\includegraphics{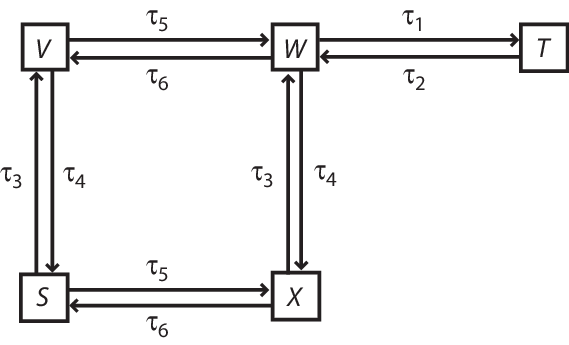}}
\caption{Digraph of a medium with set of states $\SSS=\{S,V,W,X,T\}$ and set of tokens $\TTT=\{\tau_i\}_{1\le i\le 6}$.} \label{example medium} 
\end{figure}
}

\begin{example}\label{so-basic ex} {\rm Figure \ref{example medium} displays the digraph
of a medium with set of states $\SSS =\{S,V,W,X,T\}$ and set of tokens 
$\TTT = \{\tau_i\mid 1\leq i \leq 6\}$.  It is clear that 
$\tilde\tau_1=\tau_2$, $\tilde\tau_3=\tau_4$, and $\tilde\tau_5=\tau_6$. We omit loops in digraphs representing token systems.
}
\end{example}

\begin{theorem}\label{so-ax independent} The  axioms [$\MMM$1] and [$\MMM$2] are independent.
\end{theorem}

\begin{proof}
Each of the two digraphs in Figure \ref{so-4-counter} defines a token system with the set of states $\{S,V,W\}$  satisfying one of the two axioms defining a medium. The axiom labeling each digraph indicates the failing axiom. Indeed, in Example [$\MMM$1] of Figure~\ref{so-4-counter}, there is no message producing $W$ from any other state, so Axiom [$\MMM$1] fails; Axiom [$\MMM$2] holds vacuously. In Example [$\MMM$2] of the same figure, the message $\tau_1\tau_2\tau_3$ is stepwise effective for $S$ and ineffective for $S$, but is not vacuous. Thus, Axiom [$\MMM$2] does not hold. Clearly, Axiom [$\MMM$1] holds in this case.
\end{proof}

{\begin{figure}[h!]
\centerline{\includegraphics{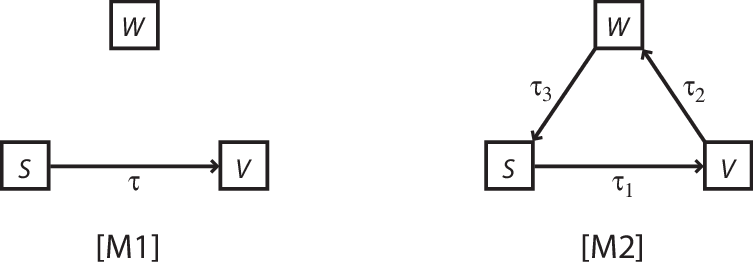}}
\caption{Digraphs of two token systems. Each digraph is labeled by the unique failing Axiom.} \label{so-4-counter} 
\end{figure}
}

\section{A `canonical' example of a medium} \label{S:example}

Let $X$ be a set and $\FFF$ be a family of subsets of $X$ such that $|\FFF|\ge 2$. For every $x\in\cup\,\FFF\setminus\cap\,\FFF$, we define transformations $\gg_x$ and $\tilde\gg_x$ of the family $\FFF$ by
$$
\gg_x:S\mapsto S\gg_x=\begin{cases}
	S\cup\{x\},	&\text{if $S\cup\{x\}\in\FFF$,}\\
	S,	&\text{otherwise,}
\end{cases}
$$
and
$$
\tilde\gg_x:S\mapsto S\tilde\gg_x=\begin{cases}
	S\setminus\{x\},	&\text{if $S\setminus\{x\}\in\FFF$,}\\
	S,	&\text{otherwise,}
\end{cases}
$$
respectively, and denote $\GGG_\FFF$ the family of all these transformations. We say that the family $\FFF$ is {\em connected} if, for any two sets $S,T\in\FFF$, there is a sequence $S_0=S,S_1,\ldots,S_n=T$ of sets in $\FFF$ such that $d(S_i,S_{i+1})=|S_i\bigtriangleup S_{i+1}|=1$ for all $i$.

\begin{lemma} \label{connected -> token system}
If $\FFF$ is connected, then $(\FFF,\GGG_\FFF)$ is a token system.
\end{lemma}

\begin{proof}
We need to show that $\gg_x\neq\tau_0$ and $\tilde\gg_x\neq\tau_0$ for any given $x$. Since $x\in\cup\,\FFF\setminus\cap\,\FFF$, there are $S,T\in\FFF$ such that $x\notin S$ and $x\in T$. Let $(S_i)$ be a sequence of sets in $\FFF$ such that $S_0=S$, $S_n=T$, and $d(S_i,S_{i+1})=1$ for all $i$. Clearly, there is $k$ such that $x\notin S_k$ and $x\in S_{k+1}$. It follows that $S_{k+1}=S_k+\{x\}$, so $S_k\gg_x=S_{k+1}$. Therefore, $\gg_x\neq\tau_0$. Evidently, $S_{k+1}\tilde\gg_x=S_k$, so $\tilde\gg_x\neq\tau_0$.
\end{proof}

\begin{remark}
{\rm The converse of the lemma does not hold. For $X=\{a,b,c,d\}$, let $\FFF=\{\{a\},\{b\},\{a,b\},\{c\},\{d\},\{c,d\}\}$. The family $\FFF$ is not connected, but $(\FFF,\GGG_\FFF)$ is a token system.
}
\end{remark}

\begin{definition}
{\rm A family $\FFF$ of subsets of a set $X$ is {\em well-graded} (a {\em wg-family}) if, for any two distinct subsets $S,T\in\FFF$ with $d(S,T)=n$, there is a sequence $S_0=S,S_1,\ldots,S_n=T$ such that $d(S_i,S_{i+1})=1$ for all $0\le i<n$.
}
\end{definition}

\begin{remark}
{\rm A family $\FFF$ of finite subsets of $X$ is well-graded if and only if the induced graph $\langle\FFF\rangle$ is an isometric subgraph of the cube $\HHH(X)$, that is, $\langle\FFF\rangle$ is a partial cube on $X$ (see~\ref{S:bipartite}).
}
\end{remark}

We will need the following result.

\begin{lemma} \label{|i-j|}
Let $(S_0,S_1,\ldots,S_n)$ be a sequence of subsets of $X$ such that 
$$
d(S_0,S_n)=n\text{~~and~~} d(S_{i-1},S_i)=1\text{~~for~~} 1\le i\le n.
$$
Then $d(S_i,S_j)=|i-j|$, for all $0\le i,j\le n$.
\end{lemma}

\begin{proof}
We may assume that $i<j$. By the triangle inequality,
\begin{align*}
n=d(S_0,S_n)&\le d(S_0,S_i)+d(S_i,S_j)+d(S_j,S_n) \\
&\le (i-0)+(j-i)+(n-j)=n.
\end{align*}
It follows that $d(S_i,S_j)=j-i$.
\end{proof}

\begin{theorem} \label{representing medium}
$(\FFF,\GGG_\FFF)$ is a medium if and only if $\FFF$ is a wg-family.
\end{theorem}

\begin{proof}
(Necessity.) Let $S$ and $T$ be two distinct sets in $\FFF$. By [$\MMM$1], there is a concise message $\bm=\tau_1\ldots\tau_n$ transforming $S$ into $T$. Let $(S_i)$ be a sequence of sets produced by $\bm$ from $S$, so $S_0=S$ and $S_n=T$. Each $\tau_i$ is either $\gg_{x_i}$ or $\tilde\gg_{x_i}$ for some $x_i$. Since $\bm$ is a concise message, all elements $x_i$ are distinct. Suppose first that $\tau_i=\gg_{x_i}$ for some $i$. Then $S_i=S_{i-1}+\{x_i\}$. Since $\bm$ is a concise message, we must have $x_i\in S_j$ for all $j\ge i$ and $x_i\notin S_j$ for all $j<i$. Hence, $x_i\in T\setminus S$. Suppose now that $\tau_i=\tilde\gg_{x_i}$ for some $i$. Then $S_i=S_{i-1}\setminus\{x_i\}$. Arguing as in the previous case, we obtain $x_i\in S\setminus T$. Therefore, $x_i\in S\bigtriangleup T$ for any $i$. On the other hand, it is clear that any element of $S\bigtriangleup T$ is one of the $x_i$'s. Thus $S\bigtriangleup T=\cup_i\{x_i\}$, so $d(S,T)=n$. Clearly, we have $d(S_{i-1},S_i)=1$, for all $i$. It follows that $\FFF$ is a wg-family.

(Sufficiency.) Let $\FFF$ be a well-graded family of subsets of some set $X$. By Lemma~\ref{connected -> token system}, $(\FFF,\GGG_\FFF)$ is a token system. It is clear that the tokens $\gg_x$ and $\tilde\gg_x$ are mutual reverses for any $x\in\cup\,\FFF\setminus\cap\,\FFF$. We need to show that Axioms [$\MMM$1] and [$\MMM$2] are satisfied for $(\FFF,\GGG_\FFF)$.

Axiom [$\MMM$1]. Let $S$ and $T$ be two distinct states in the wg-family $\FFF$, and let $(S_i)$ be a sequence of states in $\FFF$ such that $S_0=S$, $S_n=T$, $d(S,T)=n$, and $d(S_{i-1},S_i)=1$. By the last equation, for any $i$, there is $x_i$ such that $S_{i-1}\bigtriangleup S_i=\{x_i\}$. Suppose that $x_i=x_j$ for some $i<j$. We have
$$
(S_{i-1}\bigtriangleup S_j)\bigtriangleup(S_i\bigtriangleup S_{j-1})=(S_{i-1}\bigtriangleup S_i)\bigtriangleup(S_{j-1}\bigtriangleup S_j)=\{x_i\}\bigtriangleup\{x_j\}=\es.
$$
Hence, $S_{i-1}\bigtriangleup S_j=S_i\bigtriangleup S_{j-1}$, so, by Lemma~\ref{|i-j|},
$$
j-(i-1)=d(S_{i-1},S_j)=d(S_i,S_{j-1})=(j-1)-i,
$$
a contradiction. Thus, all $x_i$'s are distinct. Since $S_{i-1}\bigtriangleup S_i=\{x_i\}$, we have $S_{i-1}\tau_i=S_i$, where $\tau_i$ is either $\gg_{x_i}$ or $\tilde\gg_{x_i}$. Clearly, the message $\tau_1\ldots\tau_n$ is concise and produces $T$ from $S$.

Axiom [$\MMM$2]. Let $\bm=\tau_1\ldots\tau_n$ be a stepwise effective message for a state $S$ which is ineffective for $S$. As before, $(S_i)$ stands for the sequence of states produced by $\bm$ from $S$, so $S_0=S_n=S$. Since $S\bm=S$, for any occurrence of $\tau$ in $\bm$ there must be occurrence of $\tilde\tau$ in $\bm$. Suppose that we have two consecutive occurrences of a token $\tau=\tau_i=\tau_j=\gg_x$ in $\bm$. Then $x\in S_i$ and $x\notin S_{j-1}$. Therefore we must have an occurrence of $\tilde\tau=\tilde\gg_x$ between this two occurrences of $\tau$. A similar argument shows that there is an occurrence of a token between any two consecutive occurences of its reverse, so occurrences of token and its reverse alternate in $\bm$. Finally, let $\tau_i$ be the first occurrence of $\tau$ in $\bm$. We may assume that there are more than one occurrence of $\tau$ in $\bm$. The message $\bn=\tau_{i+1}\ldots\tau_n\tau_1\ldots\tau_i$ is stepwise effective and ineffective for $S_i$. By the previous argument, occurrences of $\tau$ and its reverse alternate in $\bn$. It follows that the number of occurrences of both $\tau$ and $\tilde\tau$ in $\bm$ is even, so $\bm$ is vacuous.
\end{proof}

Theorem~\ref{representing medium} justifies the following definition.
\begin{definition}
{\rm Let $\FFF$ be a wg-family of subsets of a set $X$. The medium $(\FFF,\GGG_\FFF)$ is said to be the {\em representing medium} of $\FFF$.
}
\end{definition}

We will show later (Theorem~\ref{isomorphism theorem}) that any medium is isomorphic to the representing medium of some wg-family of sets. This is why $(\FFF,\GGG_\FFF)$ is a `canonical' example of a medium.

The representing medium $(\BBB(X),\GGG_{\BBB(X)})$ of the family of all finite subsets of $X$ has a rather special property: 
\begin{quote}
For any state $S$ and any token $\gg$, either $\gg$ or $\tilde\gg$ is effective for $S$.
\end{quote}
Any medium satisfying this property is said to be {\em complete}.

\section{Tokens and messages of media} \label{S:tokens}

The two axioms defining a medium are quite strong. We derive a few basic consequences of these axioms. In what follows we assume that a medium $(\SSS,\TTT)$ is given.

\begin{lemma} \label{3 properties}
{\rm(i)} Any token of a medium has a reverse. In particular, if $S$ is adjacent to $V$, then $S$ and $V$ are adjacent.

{\rm(ii)} No token can be identical to its own reverse. In particular, a single token $\tau$ is a concise message for any state $S$ such that $S\tau\neq S$.

{\rm(iii)} For any two adjacent states, there is exactly one token producing one state from the other.
\end{lemma}

\begin{proof}
(i) and (ii). Let $\tau$ be a token in $\TTT$. Since $\tau\neq\tau_0$ (recall that $\tau_0$ stands for the identity transformation of $\SSS$ and is not a token), there are two distinct states $S$ and $V$ in $\SSS$ such that $S\tau=V$. By Axiom [$\MMM$1], there is a concise message $\bm$ producing $S$ from $V$. The message $\tau\bm$ is stepwise effective for $S$ and ineffective for that state. By Axiom [$\MMM$2], this message is vacuous. Hence, the message $\bm$ contains a reverse of $\tau$. It follows that there is a reverse of $\tau$ in $\TTT$. If $\tau=\tilde\tau$, then $\bm$ contains both $\tau$ and $\tilde\tau$. This contradicts the assumption that $\bm$ is a concise message.

(iii) Suppose that $S\tau_1=S\tau_2=V$, so $V$ is adjacent to $S$. By (i), the message $\tau_1\tilde\tau_2$ is stepwise effective and ineffective for $S$. By Axiom [$\MMM$2], it is vacuous, that is, $\{\tau_1,\tilde\tau_2\}$ is a pair of mutually reverse tokens. Therefore, $\tau_1=\tilde{\tilde\tau}_2=\tau_2$.
\end{proof}

Let $\tau$ be a token of a medium. We define
\begeq \label{U-set}
\UUU_\tau=\{S\in\SSS\mid S\tau\neq S\}.
\edeq
Note that $\UUU_\tau\neq\es$, since $\tau$ is a token.

\begin{lemma} \label{token-structure}
For any given $\tau\in\TTT$ we have
\roster
	\item[{\rm(i)}] $(\UUU_\tau)\tau=\UUU_{\tilde\tau}$.
	\item[{\rm(ii)}] $\UUU_\tau\cap\UUU_{\tilde\tau}=\es$.
	\item[{\rm(iii)}] The restriction $\tau|_{_{\UUU_\tau}}$ is a bijection from $\UUU_\tau$ onto $\UUU_{\tilde\tau}$ with $\tau|^{-1}_{_{\UUU_\tau}}=\tilde\tau|_{_{\UUU_{\tilde\tau}}}$.
	\item[{\rm(iv)}] $\tau$ is not a one-to-one transformation.
\endroster
\end{lemma}

\begin{proof}
(i) We have
$$
T\in(\UUU_\tau)\tau\;\eq\; S\tau=T\;(S\neq T)\;\eq\; T\tilde\tau=S\;(S\neq T)\;\eq\; T\in\UUU_{\tilde\tau}.
$$

(ii) If $S\in\UUU_\tau\cap\UUU_{\tilde\tau}$, then there exist $T\neq S$ such that $S\tau=T$ and $V\neq S$ such that $S\tilde\tau=V$, so $V\tau=S$. If $V=T$, then, by [$\MMM$2], the message $\tau\tau$ is vacuous, so $\tilde\tau=\tau$, which contradicts Lemma~\ref{3 properties}(ii). If $V\neq T$, then, by [$\MMM$1], there is a concise message $\bm$ producing $V$ from $T$. By [$\MMM$2], the message $\tau\tau\bm$ is vacuous, so we must have two occurrences of $\tilde\tau$ in $\bm$ a contradiction, since $\bm$ is a concise message. It follows that $\UUU_\tau\cap\UUU_{\tilde\tau}=\es$.

(iii) and (iv) follow immediately from (i) and (ii).
\end{proof}

\begin{lemma}
If $\bm$ is a concise message for some state $S$, then $\bm$ is effective for $S$.
\end{lemma}

\begin{proof}
If $S\bm=S$, then, by Axiom [$\MMM$2], $\bm$ must be vacuous, which contradicts our assumption that $\bm$ is a concise message.
\end{proof}

\begin{lemma}
A vacuous message $\bm$ which is stepwise effective for a state $S$ is ineffective for $S$.
\end{lemma}

\begin{proof}
Suppose that $T=S\bm\neq S$, and let $\bn$ be a concise message producing $S$ from $T$. By Axiom [$\MMM$2], the message $\bm\bn$ is vacuous, so $\bn$ must contain a pair of mutually reverse tokens, a contradiction. Hence, $S\bm=S$.
\end{proof}

%

\begin{lemma} \label{so-3-gon}
Let $S$, $V$, and $W$ be three states of the medium $(\SSS,\TTT)$ and suppose that $V=S\bm$, $W=V\bn$
for some concise messages $\bm$ and $\bn$, and $S=W\bp$ where $\bp$ is either a concise message or empty (see the diagram in Figure~\ref{triangle1}). There is at most one occurrence of each pair of mutually reverse tokens in the closed message $\bm\bn\bp$.
\end{lemma}

{\begin{figure}[h!]
\centerline{\includegraphics[scale=0.8]{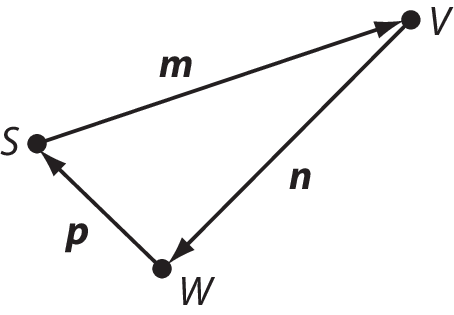}}
\caption{Diagram for Lemma~\ref{so-3-gon}.} \label{triangle1} 
\end{figure}
}

\begin{proof}
Let $\tau$ be a token in $\CCC(\bm)$. Since $\bm$ is a concise message, there is only one occurrence of $\tau$ in $\bm$ and $\tilde\tau\notin\CCC(\bm)$. By Axiom [$\MMM$2], the message $\bm\bn\bp$ is vacuous, so we must have $\tilde\tau\in\CCC(\bn)\cup\CCC(\bp)$. Suppose that $\tilde\tau\in\CCC(\bn)$ (the case when $\tilde\tau\in\CCC(\bp)\neq\es$ is treated similarly). Since $\bn$ is a concise message, there are no more occurences of $\tilde\tau$ in $\bn$ and $\tau\not\in\CCC(\bn)$. Thus there is only one occurrence of the pair $\{\tau,\tilde\tau\}$ in the message $\bm\bn$. The pair $\{\tau,\tilde\tau\}$ cannot occur in $\bp$, since $\bp$ is a concise message. The result follows.
\end{proof}

\begin{corollary} \label{2-gon}
Let $\bm$ and $\bn$ be two concise messages producing $V$ from $S$. Then the string $\bn$ is a permutation of the string $\bm$. In particular, $\ell(\bm)=\ell(\bn)$.
\end{corollary}

One can say more in the special case when $\bp=\tau$ is a single token.

\begin{lemma} \label{small triangle}
Let $S$, $V$ and $W$ be distinct states of a medium and suppose that
$$
V=S\bm,\quad W=V\bn,\quad S=W\tau
$$
for some concise messages $\bm$ and $\bn$ and a token $\tau$ (see Figure~\ref{triangle2}). Then 
$$
\tau\notin\CCC(\bn),\quad \tau\notin\CCC(\bm),
$$
and either
$$
\tilde\tau\in\CCC(\bm),\quad\text{$\bn\tau$ is a concise message},\quad \CCC(\bn\tau)=\CCC(\widetilde\bm),\quad \ell(\bm)=\ell(\bn)+1,
$$
or
$$
\tilde\tau\in\CCC(\bn),\quad\text{$\tau\bm$ is a concise message},\quad \CCC(\tau\bm)=\CCC(\widetilde\bn),\quad \ell(\bn)=\ell(\bm)+1.
$$
Accordingly,
\begeq \label{m-n=1}
|\ell(\bm)-\ell(\bn)|=1.
\edeq
\end{lemma}

{\begin{figure}[h!]
\centerline{\includegraphics[scale=0.8]{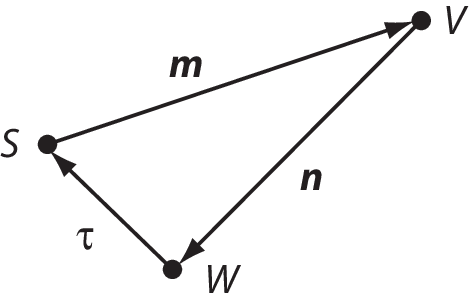}}
\caption{For Lemma~\ref{small triangle}.} \label{triangle2} 
\end{figure}
}

\begin{proof}
By Lemma~\ref{so-3-gon}, $\tau\notin\CCC(\bn)$, $\tau\notin\CCC(\bm)$, and $\tilde\tau$ occurs either in $\bm$ or in $\bn$. Suppose that $\tilde\tau\in\CCC(\bm)$. By the same lemma, neither $\tau$ not $\tilde\tau$ occurs in $\bn$. Therefore, $\bn\tau$ is a concise message. The equality $\CCC(\bn\tau)=\CCC(\widetilde\bm)$ also follows from Lemma~\ref{so-3-gon}. Since $\bm$ is a concise message, we have
$$
\ell(\bm)=|\CCC(\bm)|=|\CCC(\widetilde\bm)|=|\CCC(\bn\tau)|=\ell(\bn)+1.
$$
The case when $\tilde\tau\in\CCC(\bn)$ is treated similarly.
\end{proof}

\begin{remark}
{\rm In each of two mutually exclusive cases of Lemma~\ref{small triangle} there are $2$-gons constructed from concise messages $\bm$, $\bn$, and $\tau$. For instance,

if $\tilde\tau\in\CCC(\bm)$, then $\bm\bn\tau$ is a $2$-gon for $S$,

if $\tilde\tau\in\CCC(\bn)$, then $\tau\bm\bn$ is a $2$-gon for $W$.
}
\end{remark}

\begin{remark}
{\rm Equation~(\ref{m-n=1}) also follows from Lemma~\ref{tr-main}.
}
\end{remark}
The results of Lemma~\ref{so-3-gon} suggest an interpretation of the length function on messages. First, by Corollary~\ref{2-gon}, we have $\ell(\bm)=\ell(\bn)$ for any two concise messages $\bm$ and $\bn$ producing a state $V$ from a state $S$. Therefore the function
\begeq \label{d-function}
\dd(S,V)=\begin{cases}
	\ell(\bm), &\text{if $S\bm=V$,}\\
	0, &\text{otherwise,}
\end{cases}
\edeq
where $\bm$ is a concise message, is well-defined. Since $\ell(\widetilde\bm)=\ell(\bm)$, the function $\dd$ is symmetric.

Second, for the messages in Figure~\ref{triangle1} we have
$$
|\CCC(\bp)|\le|\CCC(\bn)|+|\CCC(\bm)|,
$$
by Lemma~\ref{so-3-gon}. Indeed, for any token in $\CCC(\bp)$ we have a unique matching reverse either in $\CCC(\bm)$ or in $\CCC(\bn)$. Since the length of a concise message equals the cardinality of its content, we have the triangle inequality
$$
\dd(S,W)\le\dd(S,V)+\dd(V,W).
$$
It is easy to verify that this inequality holds for any choice of states $S,V,W\in\SSS$.

We obtained the following result.

\begin{theorem} \label{media metric}
The function $\dd(S,V)$ defined by~{\rm(\ref{d-function})} is a metric on $\SSS$.
\end{theorem}

\section{Bipartite graphs and partial cubes} \label{S:bipartite}

Some useful properties of media can be derived from metric properties of their graphs. First, we formulate two well-known characterization properties of bipartite graphs (see, for instance,~\cite{aA98}). In what follows, $\dd$ stands for the graph distance.

\begin{theorem}
A graph $G$ is bipartite if and only if it contains no closed walk of odd length.
\end{theorem}

\begin{theorem} \label{bipartite triangle}
A connected graph is bipartite if and only if for every vertex $T$ there is no edge $\{S,V\}$ such that $\dd(T,S)=\dd(T,V)$.
\end{theorem}

\begin{lemma} \label{tr-main}
Let $\{S,V\}$ be an edge of a connected bipartite graph $G$ and $W$ be a vertex of $G$. Then
\begeq
|\dd(W,S)-\dd(W,V)|=1.
\edeq
\end{lemma}

\begin{proof}
By the triangle inequality,
$$
|\dd(W,S)-\dd(W,V)|\le 1.
$$
By Theorem~\ref{bipartite triangle}, $\dd(W,S)\neq\dd(W,V)$. Since values of the function $\dd$ are whole numbers, we have~(\ref{tr-main}).
\end{proof}

{\begin{figure}[h!]
\centerline{\includegraphics{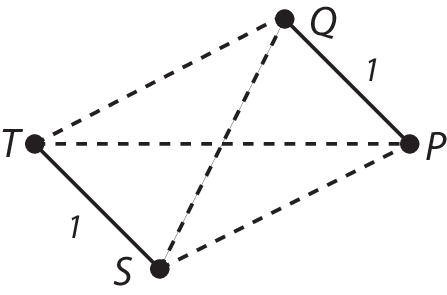}}
\caption{Two edges of a graph. Note that $\dd(S,T)=\dd(P,Q)=1$.} \label{tetrahedron} 
\end{figure}
}

\begin{theorem} \label{tetrahedron theorem}
Let $G$ be a connected bipartite graph and $\{S,T\}$ and $\{P,Q\}$ be two distinct edges of $G$(see Figure~\ref{tetrahedron}). There are six mutually exclusive, exhaustive cases:

Case 1: $\quad\dd(T,P)=\dd(S,Q)=\dd(S,P)+1=\dd(T,Q)-1$

Case 2: $\quad\dd(T,P)=\dd(S,Q)=\dd(S,P)-1=\dd(T,Q)+1$

Case 3: $\quad\dd(S,P)=\dd(T,Q)=\dd(T,P)+1=\dd(S,Q)-1$

Case 4: $\quad\dd(S,P)=\dd(T,Q)=\dd(T,P)-1=\dd(S,Q)+1$

Case 5: $\quad\dd(S,P)=\dd(T,Q)=\dd(T,P)+1=\dd(S,Q)+1$

Case 6: $\quad\dd(S,P)=\dd(T,Q)=\dd(T,P)-1=\dd(S,Q)-1$
\end{theorem}

\begin{proof}
By applying the result of Lemma~\ref{tr-main} to `triangles' $PST$, $PQT$, $PSQ$, and $STQ$, we obtain four equations:
\begin{align}
|\dd(S,P)-\dd(T,P)|&=1 \label{tr-PST} \\
|\dd(T,Q)-\dd(T,P)|&=1 \label{tr-PQT} \\
|\dd(S,P)-\dd(S,Q)|&=1 \label{tr-PSQ} \\
|\dd(T,Q)-\dd(S,Q)|&=1 \label{tr-STQ}, 
\end{align}
It is helpful to regard the absolute values in these equations as distances on the number line.

Suppose that $\dd(S,P)\neq\dd(T,Q)$. Then, $2\dd(T,P)=\dd(S,P)+\dd(T,Q)$, by~(\ref{tr-PST}) and~(\ref{tr-PQT}), and $2\dd(S,Q)=\dd(S,P)+\dd(T,Q)$, by~(\ref{tr-PSQ}) and~(\ref{tr-STQ}). Therefore, $\dd(T,P)=\dd(S,Q)$. There are two mutually exclusive possibilities in this case:

Case 1: $\quad\dd(T,P)=\dd(S,Q)=\dd(S,P)+1=\dd(T,Q)-1$

Case 2: $\quad\dd(T,P)=\dd(S,Q)=\dd(S,P)-1=\dd(T,Q)+1$

Suppose that $\dd(S,P)=\dd(T,Q)$ and $\dd(T,P)\neq\dd(S,Q)$. The same argument as above, applied to equations~(\ref{tr-PST}),~(\ref{tr-PSQ}) and~(\ref{tr-PQT}),~(\ref{tr-STQ}), shows that there are again two mutually exclusive possibilities:

Case 3: $\quad\dd(S,P)=\dd(T,Q)=\dd(T,P)+1=\dd(S,Q)-1$

Case 4: $\quad\dd(S,P)=\dd(T,Q)=\dd(T,P)-1=\dd(S,Q)+1$

Finally, suppose that $\dd(S,P)=\dd(T,Q)$ and $\dd(T,P)=\dd(S,Q)$. Obviously, we have either

Case 5: 
\begeq \label{case5}
\quad\dd(S,P)=\dd(T,Q)=\dd(T,P)+1=\dd(S,Q)+1
\edeq

\noindent
or

Case 6: 
\begeq \label{case6}
\quad\dd(S,P)=\dd(T,Q)=\dd(T,P)-1=\dd(S,Q)-1
\edeq

It is clear that the six cases are mutually exclusive and exhaustive.
\end{proof}

In the first four cases we obtain the identities:

Case 1: $\quad\dd(T,Q)=\dd(S,P)+2$.

Case 2: $\quad\dd(S,P)=\dd(T,Q)+2$.

Case 3: $\quad\dd(S,Q)=\dd(T,P)+2$.

Case 4: $\quad\dd(T,P)=\dd(S,Q)+2$.

It follows that in these four cases the edges $\{S,T\}$ and $\{P,Q\}$ belong to a shortest path in $G$ (see Figure~\ref{cases1-4}).

The remaining two cases are depicted in Figure~\ref{cases5-6}. 
In these two cases, the four vertices do not belong to any shortest path in $G$. It is natural to call the configuration defined by these vertices a `rectangle'---the opposite `sides' are equal as well as the two `diagonals'.

{\begin{figure}[h!]
\centerline{\includegraphics{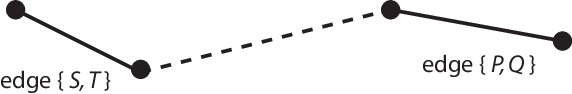}}
\caption{A path for Cases 1--4 containing edges $\{S,T\}$ and $\{P,Q\}$.} \label{cases1-4} 
\end{figure}
}

{\begin{figure}[h!]
\centerline{\includegraphics{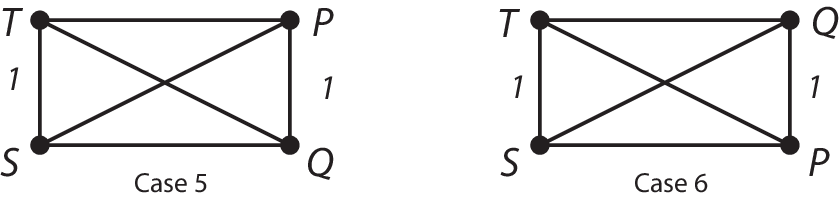}}
\caption{Two `rectangles' for Cases 5 and 6.} \label{cases5-6} 
\end{figure}
}

It is easy to verify that the first four cases can be distinguished from the last two by the following conditions:
\begin{align*}
\dd(T,P)+\dd(S,Q)&=\dd(S,P)+\dd((T,Q)\quad\text{in Cases 1--4,}\\
\dd(T,P)+\dd(S,Q)&\neq\dd(S,P)+\dd((T,Q)\quad\text{in Cases 5 and 6.}
\end{align*}

Let us recall~\cite{wI00,pW84} that{\em Winkler's relation $\Theta$} is a binary relation on the set of edges $E$ defined by
$$
\{S,T\}\Theta\{P,Q\}\quad\EQ\quad\dd(S,P)+\dd(T,Q)\neq\dd(S,Q)+\dd(T,P).
$$

We established, in particular, the following result (cf.~\cite[Lemma~2.2]{wI00}).

\begin{theorem}
Two edges of a connected bipartite graph stand in Winkler's relation $\Theta$ if and only if they do not belong to the same shortest path.
\end{theorem}


\begin{definition}
{\rm A {\em cube on a set $X$}, $\HHH(X)$, has the set $\BBB(X)$ of all finite subsets of $X$ as the set of vertices; $\{S,T\}$ is an edge of $\HHH(X)$ if $|S\bigtriangleup T|=1$.
A {\em partial cube} is a graph that is isometrically embeddable into some cube $\HHH(X)$.
}
\end{definition}

\begin{definition} \label{semicube graph}
{\rm Let $G=(V,E)$ be a connected graph and $\dd$ be the graph distance on $G$. For any $\{S,T\}\in E$, the sets
$$
W_{ST}=\{P\in V\mid \dd(V,S)<\dd(V,T)\}
$$
are called {\em semicubes} of $G$. The semicubes $W_{ST}$ and $W_{TS}$ are called {\em opposite semicubes of $G$}.
}
\end{definition}

It is easy to prove (see~\cite{sO06a}) that a graph $G$ is bipartite if and only if the opposite semicubes $W_{ST}$ and $W_{TS}$ form a partition of $V$.

The following theorem (Theorem~2.10 in~\cite{wI00}) summarizes two major characterizations of partial cubes due to Djokovi\'{c}~\cite{dD73} and Winkler~\cite{pW84}.

\begin{theorem} \label{DWT}
For a connected graph $G=(V,E)$ the following statements are equivalent:
\roster
	\item[{\rm(i)}] $G$ is a partial cube.
	\item[{\rm(ii)}] $G$ is bipartite and for every edge $\{S,T\}$ the semicube $W_{ST}$ induces a convex subgraph of $G$.
	\item[{\rm(iii)}] $G$ is bipartite and $\Theta$ is an equivalence relation on $E$.
\endroster
\end{theorem}

We give another useful characterization of partial cubes.

\begin{theorem} \label{labels}
A graph $G=(V,E)$ is a partial cube if and only if it is possible to label its edges by elements of some set $J$ such that
\roster
	\item[\rm(i)] Edges of any shortest path of $G$ are of different labels.
	\item[\rm(ii)] In each closed walk of $G$ every label appears an even number of times.
\endroster
\end{theorem}

\begin{proof}
(Necessity.) Without loss of generality, we may assume that $G=(\FFF,\EEE)$ is an isometric subgraph of a cube $\HHH(J)$ such that $\cap\,\FFF=\es$ and $\cup\;\FFF=J$ for a wg-family $\FFF$. For any edge $\{S,T\}$ of $G$ there is an element $j\in J$ such that $S\bigtriangleup T=\{j\}$, so we can label edges of $G$ by elements of $J$.

(i) Let $S_0=S,S_1,\ldots,S_n=T$ be a shortest path from $S$ to $T$ in $G$. For every $i$, we have $S\cap T\SB S_i\SB S\cup T$. Therefore,
$$
\{j_i\}=S_{i-1}\bigtriangleup S_i\SB S\bigtriangleup T.
$$
Since $(S_i)$ is a shortest path, $|S\bigtriangleup T|=d(S,T)=n$. It follows that all labels $j_i$ are distinct.

(ii) Let $S_0,S_1,\ldots,S_n=S_0$ be a closed walk $W$ in $G$ and let $E_p=\{S_{p-1},S_p\}$ be the first edge in $W$ labeled by $j$, so $S_{p-1}\bigtriangleup S_p=\{j\}$. We assume that $j\notin S_{p-1}$ and $j\in S_p$; the other case is treated similarly. Since $E_p$ is the first edge of $W$ labeled by $j$, we must have $j\notin S_0$. Since the walk $W$ is closed and $j\in S_p$, we must have another occurrence of $j$ in $W$. Let $E_q=\{S_{q-1},S_q\}$ be the next edge of $W$ labeled by $j$. We have $j\in S_{q-1}$ and $j\notin S_q$. By repeating this argument, we partition the occurrences of $j$ in $W$ into pairs, so the total number of these occurrences must be even.

(Sufficiency.) Let $S_0$ be a fixed vertex of $G$. For any vertex $S\in V$ and a shortest path $p$ from $S_0$ to $S$, we define
$$
J_S=\{j\in J\mid \text{$j$ is a label of an edge of $p$}\},
$$
and $J_{S_0}=\es$. The set $J_S$ is well-defined. Indeed, let $q$ be another shortest path from $S_0$ to $S$ and $\tilde q$ be its reverse, so $p\tilde q$ is a closed walk. By (i) and (ii), $J_S$ does not depend on the choice of $p$.

The correspondence $\aa:S\mapsto J_S$ defines an isometric embedding of $G$ into the cube $\HHH(J)$. Indeed, for $S,T\in V$, let $p$ (resp. $q$) be a shortest path from $S_0$ to $S$ (resp. $T$) and let $r$ be a shortest path from $S$ to $T$. By (ii) applied to the closed walk $pr\tilde q$ and (i), we have
$$
j\in J_S\bigtriangleup J_T\quad\EQ\quad\text{$j$ is a label of an edge of $r$},
$$
so $\dd(S,T)=|J_S\bigtriangleup J_T|=d(J_S,J_T)$.
\end{proof}

\section{The graph of a medium} \label{S:graphs}

\begin{definition}
{\rm The {\em graph} of a medium $(\SSS,\TTT)$ has $\SSS$ as the set of vertices; two vertices of the graph are adjacent if and only if the corresponding states are adjacent in the medium.
}
\end{definition}

By Lemma~\ref{3 properties}, for any two adjacent states $S$ and $T$ of a medium $(\SSS,\TTT)$ there is a unique token $\tau$ such that $S\tau=T$ and $T\tilde\tau=S$. Thus, a unique pair of mutually reversed tokens $\{\tau,\tilde\tau\}$ is assigned to each edge $\{S,V\}$ of the graph of $(\SSS,\TTT)$.

Let $(\SSS,\TTT)$ be a medium and $G$ be its graph. If $\bm=\tau_1\ldots\tau_m$ is a stepwise effective message for a state $S$ producing a state $V$, then the sequence of vertices $(S_i)$ of $G$ produced by $\bm$, is a walk in $G$; the vertex $S_0=S$ is a tail of this walk and the vertex $S_m=V$ is its head. On the other hand, if the sequence of vertices $S_0=S,S_1,\ldots,S_m=V$ is a walk in $G$, then edges $\{S_{i-1},S_i\}$ define unique tokens $\tau_i$ such that $S_{i-1}\tau_i=S_i$. Then $\bm=\tau_1\ldots\tau_m$ is a stepwise effective message for the state $S$ producing the state $V$. Thus we have a one-to-one correspondence between stepwise effective messages of the medium and walks in its graph. In particular, a closed message for some state produces a closed walk in $G$.

A deeper connection between media and their graphs is the result of the following theorem.

\begin{theorem} \label{concise=shortest}
Let $(\SSS,\TTT)$ be a medium and $G$ be its graph. If $\bm=\tau_1\ldots\tau_m$ is a concise message producing a state $V$ from a state $S$, then the sequence of vertices $(S_i)$ produced by $\bm$ forms a shortest path connecting $S$ and $V$ in the graph $G$. Conversely, if $S_0=S,S_1,\ldots,S_m=V$ is a shortest path in $G$, then the corresponding message is a concise message of $G$.
\end{theorem}

\begin{proof}
(Necessity.) Let $P_0=S,P_1,\ldots,P_n=V$ be a path in $G$ joining $S$ to $V$ and $\bn=\mu_1\ldots\mu_n$ be the (stepwise effective) message of the medium corresponding to this path. By Axiom [$\MMM$2], the message $\bm\widetilde\bn$ is vacuous, so $\ell(\bm)\le\ell(\widetilde\bn)=\ell(\bn)$, since $\bm$ is a concise message for $S$. Thus the sequence $(S_i)$ is a shortest path in $G$.

(Sufficiency.) Let $S_0=S,S_1,\ldots,S_m=V$ be a shortest path in $G$ and let $\bm=\tau_1\ldots\tau_m$ be the corresponding stepwise effective message of the medium. By Axiom [$\MMM$1], there is a concise message $\bn$ producing $V$ from $S$. By the necessity part of the proof, the walk defined by $\bn$ is a shortest path from $S$ to $V$, so $\ell(\bn)=\ell(\bm)$. By Axiom [$\MMM$2], the message $\bm\tilde\bn$ must be vacuous. Since the message $\bn$ is concise and $\ell(\bn)=\ell(\bm)$, the message $\bm$ must be concise.
\end{proof}

Let $G$ be the graph of a medium. By Axiom [$\MMM$1], $G$ is connected. Let $S_0=S,S_1,\ldots,S_n=S$ be a closed walk. By Axiom [$\MMM$2], the corresponding message of the medium is vacuous. Therefore it must be of even length. It follows (see~\ref{S:bipartite}) that the graph of a medium is bipartite. Note that not every connected bipartite graph is the graph of some medium. 

\begin{example} \label{bipartite not medium}
{\rm The simplest counterexample is the complete bipartite graph $K_{2,3}$ shown in Figure~\ref{K-23}. Suppose that this graph is the graph of a medium and let $\tau$ be a token producing $T$ from $S$. By Axiom [$\MMM$2], the closed message producing the sequence of states $(S,T,P,V,S)$ must be vacuous and therefore contain an occurrence of $\tilde\tau$. We cannot have $T\tilde\tau=P$ or $V\tilde\tau=S$, since tokens are functions. Therefore, $P\tilde\tau=V$, so $V\tau=P$. The same argument applied to the closed message producing the sequence $(V,P,Q,S,V)$ shows that $S\tau=Q$. Thus $S\tau=T$ and $S\tau=Q$, a contradiction.
}
\end{example}

{\begin{figure}[h!]
\centerline{\includegraphics{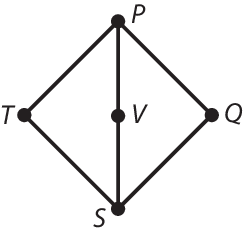}}
\caption{Complete bipartite graph $K_{2,3}$.} \label{K-23} 
\end{figure}
}

It follows from Theorem~\ref{concise=shortest} that the metric $\dd$ on the set of states of a medium is the graph distance on the graph of that medium.

\section{Contents} \label{S:contents}

\begin{definition}
{\rm Let $(\SSS,\TTT)$ be a medium. For any state $S$, the {\em content} of $S$ is the set $\widehat{S}$ of all tokens each of which is contained in at least one concise message producing $S$. The family
$\widehat{\SSS}=\{\widehat{S}\,|\,S\in\SSS\}$ is called the \emph{content family} of $\SSS$.
}
\end{definition}

\begin{lemma}
The content of a state cannot contain both a token and its reverse.
\end{lemma}

\begin{proof}
Suppose that $S\bm=W\bn=V$ for two concise messages $\bm$ and $\bn$ and let $\bp$ be a concise message producing $S$ from $W$, if $W\neq S$, and empty, if $W=S$. By Lemma~\ref{so-3-gon}, there is at most one occurrence of any token $\tau$ in the message $\bm\widetilde\bn\bp$. Therefore we cannot have both $\tau\in\CCC(\bm)$ and $\tilde\tau\in\CCC(\bn)$.
\end{proof}

\begin{theorem} \label{Theorem1.17-1}
For any token $\tau$ and any state $S$, we have either $\tau\in\widehat{S}$ or $\tilde{\tau}\in\widehat{S}$. Consequently, $|\widehat{S}|=|\widehat{V}|$ for any two states $S$ and $V$.
\end{theorem}

\begin{proof} 
Since $\tau$ is a token, there are two states $V$ and $W$ such that $W=V\tau$. By Axiom [M1], there are concise messages $\bm$ and $\bn$ such that $S=V\bm$ and $S=W\bn$. By Lemma~\ref{small triangle}, there are two mutually exclusive options: either $\tilde\tau\in\CCC(\bn)$ or $\tau\in\CCC(\bm)$.
\end{proof}

\begin{theorem} \label{Theorem1.16}
If $S$ and $V$ are two distinct states, with $S\bm=V$ for some concise message $\bm$, then \mbox{$\widehat{V}\setminus\widehat{S}=\CCC(\bm)$}.
\end{theorem}

\begin{proof} 
Let $\tau$ be a token in $\CCC(\bm)$, so $\tilde\tau\in\CCC(\widetilde\bm)$. Thus, $\tau\in\widehat{V}$ and $\tilde\tau\in\widehat{S}$. By Theorem~\ref{Theorem1.17-1}, $\tau\notin\widehat{S}$. It follows that $\tau\in\widehat{V}\setminus\widehat{S}$, that is, $\CCC(\bm)\SB\widehat{V}\setminus\widehat{S}$.

If $\tau\in\widehat{V}\setminus\widehat{S}$, then $\tau\in\widehat{V}$ and $\tau\notin\widehat{S}$, so, by Theorem~\ref{Theorem1.17-1}, $\tilde\tau\in\widehat{S}$. Since $\tau\in\widehat{V}$, there is a concise message $\bn$ producing the state $V$ from some state $W$ such that $\tau\in\CCC(\bn)$, so $\tilde\tau\in\CCC(\widetilde\bn)$. Let $\bp$ be a concise message producing $S$ from $W$ (or empty if $S=W$). By Lemma~\ref{so-3-gon}, there is exactly one occurrence of the pair $\{\tau,\tilde\tau\}$ in the message $\bm\widetilde\bn\bp$. Since $\tilde\tau\in\widehat{S}$, we have $\tau\notin\CCC(\bp)$. Hence, $\tau\in\CCC(\bm)$. In both cases we have $\widehat{V}\setminus\widehat{S}\SB\CCC(\bm)$. The result follows.
\end{proof}

\begin{theorem} \label{Theorem1.17-2}
For any two states $S$ and $V$ we have 
$$
S=V\quad\eq\quad\widehat{S}=\widehat{V}.
$$
\end{theorem}

\begin{proof}
Suppose that $\widehat{S}=\widehat{V}$, $S\neq V$, and let $\bm$ be a concise message producing $V$ from $S$. By Theorem~\ref{Theorem1.16},
$$
\es=\widehat{V}\setminus\widehat{S}=\CCC(\bm),
$$
a contradiction. Thus, $\widehat{S}=\widehat{V}\;\imp\;S=V$. The implication $S=V\;\imp\;\widehat{S}=\widehat{V}$ is trivial.
\end{proof}

\begin{theorem} \label{Theorem1.14}
Let $\bm$ and $\bn$ be two concise messages transforming some state $S$. Then $S\bm=S\bn$ if and only if $\CCC(\bm)=\CCC(\bn)$.
\end{theorem}

\begin{proof}
(Necessity.) Suppose that $V=S\bm=S\bn$. By Theorem~\ref{Theorem1.16},
$$
\CCC(\bm)=\widehat{V}\setminus\widehat{S}=\CCC(\bn).
$$

(Sufficiency.) Suppose that $\CCC(\bm)=\CCC(\bn)$ and let $V=S\bm$ and $W=S\bn$. By Theorem~\ref{Theorem1.16},
$$
\widehat{V}\Delta\widehat{S}=\CCC(\bm)\cup\CCC(\widetilde{\bm})=\CCC(\bn)\cup\CCC(\widetilde{\bn})=\widehat{W}\Delta\widehat{S},
$$
which implies $\widehat{V}=\widehat{W}$. By Theorem~\ref{Theorem1.17-2}, $V=W$.
\end{proof}

\begin{definition} \label{semicube medium}
{\rm Let $\tau$ be a token of a medium $(\SSS,\TTT)$. The subset
$$
\WWW_\tau=\{S\in\SSS\mid\tau\in\widehat{S}\}
$$
of $\SSS$ is called a {\em semicube} of the medium $(\SSS,\TTT)$. Semicubes $\WWW_\tau$ and $\WWW_{\tilde\tau}$ are called {\em opposite semicubes}.
}
\end{definition}

By Theorem~\ref{Theorem1.17-1}, we have
$$
\WWW_\tau\cap\WWW_{\tilde\tau}=\es\quad\text{and}\quad\WWW_\tau\cup\WWW_{\tilde\tau}=\SSS,
$$
so opposite cubes form a bipartition of the set $\SSS$. It is also clear that
$$
\UUU_\tau\SB\WWW_{\tilde\tau}\qquad\text{for any $\tau\in\TTT$,}
$$
where $\UUU_\tau$ is defined by equation~(\ref{U-set}).

\begin{lemma}
Let $S$ and $T$ be two distinct states in $\WWW_\tau$ and $\bm=\tau_1\ldots\tau_n$ be a concise message transforming $S$ into $T$. All states in the sequence $(S_i)$ produced by $\bm$ from $S$ belong to the set $\WWW_\tau$.
\end{lemma}

\begin{proof}
Suppose to the contrary that there are states in $(S_i)$ that belong to the semicube $\WWW_{\tilde\tau}$. Let $j$ be the first index such that $S_j\in\WWW_\tau$ and $S_{j+1}\in\WWW_{\tilde\tau}$, so $\tau\in\widehat{S}_j$ and $\tilde\tau\in\widehat{S}_{j+1}$. By Theorem~\ref{Theorem1.17-1}, $\tilde\tau\notin\widehat{S}_j$, and, by Theorem~\ref{Theorem1.16}, $\tau_{j+1}=\tilde\tau$, since $S_j\tau_{j+1}=S_{j+1}$. Let now $k$ be the first index such that $S_k\in\WWW_{\tilde\tau}$ and $S_{k+1}\in\WWW_\tau$. By repeating the previous argument, we obtain $\tau_{k+1}=\tau$, which contradicts our assumption that $\bm$ is a concise message.
\end{proof}

By Theorem~\ref{concise=shortest} and previous lemma, we have the following result.

\begin{theorem} \label{convex semicubes}
Let $G$ be the graph of a medium $(\SSS,\TTT)$. For any token $\tau\in\TTT$, the subgraph induced by the semicube $\WWW_\tau$ is convex.
\end{theorem}

The semicubes of a medium can be metrically characterized as follows.

\begin{theorem} \label{metric semicubes}
Let $(\SSS,\TTT)$ be a medium. For any $\tau\in\TTT$ and $S,T\in\SSS$ such that $S\tau=T$,
$$
\WWW_\tau=\{V\in\SSS \mid \dd(V,T)<\dd(V,S)\}.
$$
\end{theorem}

\begin{proof}
Let $\bm$ and $\bn$ be concise messages producing $S$ and $T$, respectively, from $V\in\WWW_{\tau}$. We have $\tau\in\widehat{V}$, which implies $\tilde\tau\notin\widehat{V}$. Similarly, $T\in\WWW_{\tau}$ implies $\tau\in\widehat{T}$ and $\tilde\tau\notin\widehat{T}$. Therefore, $\tau,\tilde\tau\notin\widehat{T}\setminus\widehat{V}=\CCC(\bn)$. It follows that $\bn\tilde\tau$ is a concise message producing $S$ from $V$. By Theorem~\ref{Theorem1.14}, $\CCC(\bn\tilde\tau)=\CCC(\bm)$. It follows that $\dd(V,T)=\ell(\bm)<\ell(\bn)=\dd(V,S)$. A similar argument shows that $\dd(V,S)<\dd(V,T)$, if $V\in\WWW_{\tilde\tau}$. The result follows.
\end{proof}

If $G$ is the graph of a medium, then $\WWW_\tau=W_{ST}$ for any two states $S$ and $T$ such that $S\tau=T$. ($W_{ST}$ is the semicube of the graph $G$; see~\ref{S:bipartite}.) Clearly, in this case, $\WWW_{\tilde\tau}=W_{TS}$. By Theorem~\ref{convex semicubes}, the semicubes of $G$ are convex.

\section{Closed messages of media} \label{S:circuits}

The structure of a $2$-gon $\bm\bn$ for a state $S$ is determined by the result of Corollary~\ref{2-gon}---the string $\widetilde\bn$ is a permutation of the string $\bm$ and therefore the concise messages $\bm$ and $\widetilde\bn$ have the same content. Note that the message $\bn\bm$ is a $2$-gon for the state $V=S\bm$. 

In this section, we are concerned with closed messages of a medium $(\SSS,\TTT)$ that can be constructed by using a small number of concise messages. A $2$-gon is an example of such a closed message. Another example is the closed message $\bm\bn\bp$ ($3$-gon) in Figure~\ref{triangle1}. The structure of this closed message is described in Lemma~\ref{so-3-gon}. 

{\begin{figure}[h!]
\centerline{\includegraphics{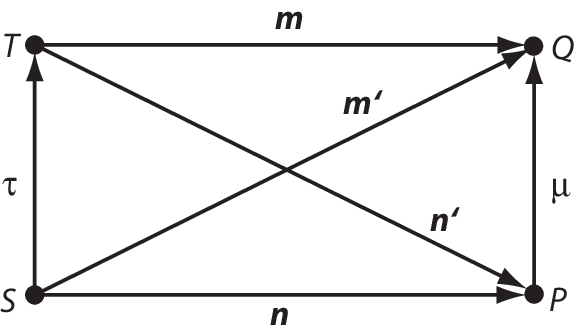}}
\caption{Four states and six concise messages.} \label{rectangle2} 
\end{figure}
}

We apply the results of Theorem~\ref{tetrahedron theorem} in~\ref{S:bipartite} to media and consider four distinct states $S$, $T$, $P$, and $Q$ and six concise messages $\tau$, $\mu$, $\bm$, $\bm'$, $\bn$, and $\bn'$ such that
$$
S\tau=T,\quad P\mu=Q,\quad T\bm=S\bm'=Q,\quad\text{and}\quad S\bn=T\bn'=P
$$
(see Figure~\ref{rectangle2}).

Since the graph of a medium is connected and bipartite, we have the following result. (In the proof we use the results of Theorem~\ref{concise=shortest}.)

\begin{theorem} \label{6 cases}
There are six mutually exclusive, exhaustive cases for the diagram in Figure~\ref{rectangle2}:

Case 1: $\quad\tilde\tau\bn\mu$ is a concise message for $T$ with $\CCC(\tilde\tau\bn\mu)=\CCC(\bm)$.

Case 2: $\quad\tau\bm\tilde\mu$ is a concise message for $S$ with $\CCC(\tau\bm\tilde\mu)=\CCC(\bn)$.

Case 3: $\quad\tau\bn'\mu$ is a concise message for $S$ with $\CCC(\tau\bn'\mu)=\CCC(\bm')$.

Case 4: $\quad\tilde\tau\bm'\tilde\mu$ is a concise message for $T$ with $\CCC(\tilde\tau\bm'\tilde\mu)=\CCC(\bn')$.

\noindent
Accordingly, $\tau\neq\mu$, $\tau\neq\tilde\mu$, $\CCC(\bm)\neq\CCC(\bn)$, and $\ell(\bm)+\ell(\bn)=\ell(\bm')+\ell(\bn')$ in all four cases listed above.

Case 5: $\quad\tau\bn'$, $\bn'\mu$, $\tilde\tau\bm'$, $\bm'\tilde\mu$ are concise messages satisfying conditions
$$
\CCC(\bm)=\CCC(\bn'\mu)=\CCC(\tilde\tau\bm')\quad\text{and}\quad\CCC(\bn)=\CCC(\tau\bn')=\CCC(\bm'\tilde\mu).
$$
Accordingly, $\tau=\tilde\mu$, $\CCC(\bm)=\CCC(\bn)$, and $\ell(\bm)+\ell(\bn)=\ell(\bm')+\ell(\bn')+2$.

Case 6: $\quad\tau\bm$, $\bn\mu$, $\tilde\tau\bn$, $\bm\tilde\mu$ are concise messages satisfying conditions
\begeq \label{eq6}
\CCC(\bm')=\CCC(\tau\bm)=\CCC(\bn\mu)\quad\text{and}\quad\CCC(\bn')=\CCC(\bm\tilde\mu)=\CCC(\tilde\tau\bn).
\edeq
Accordingly, 
\begeq \label{eqs6}
\tau=\mu,\quad\CCC(\bm)=\CCC(\bn),\quad\text{and}\quad\ell(\bm)+\ell(\bn)+2=\ell(\bm')+\ell(\bn').
\edeq
\end{theorem}

\begin{proof}
Clearly, the six cases of the theorem are mutually exclusive and exhaustive. The proofs for the first four cases are straightforward and omitted. The proofs of the two remaining cases are very similar, so we prove only the last case.

Since $\dd(S,Q)=\dd(T,Q)+1$, the path corresponding to the stepwise effective message $\tau\bm$ is a shortest path. It follows that $\tau\bm$ is a concise message producing $Q$ from $S$. Since $S\bm'=Q$, we have $\CCC(\bm')=\CCC(\tau\bm)$. Similar arguments show that $\bn\mu$, $\tilde\tau\bn$, and $\bm\tilde\mu$ are concise messages satisfying equations~(\ref{eq6}).

Suppose that $\tau\neq\mu$. By the first equation in~(\ref{eq6}), $\tau\in\CCC(\bn)$, which is impossible because $\tilde\tau\bn$ is a concise message. Hence, $\tau=\mu$. The remaining equations in~(\ref{eqs6}) follow immediately from~(\ref{eq6}).
\end{proof}

{\begin{figure}[h!]
\centerline{\includegraphics{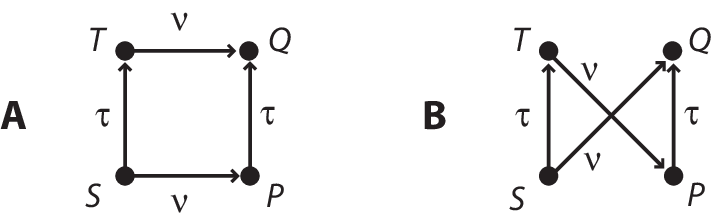}}
\caption{{\bf A} Case 6. {\bf B} Case 5, ``twisted'' Case 6.} \label{2examples} 
\end{figure}
}

\begin{example}
{\rm Cases 6 and 5 are illustrated by drawings {\bf A} and {\bf B}, respectively, in Figure~\ref{2examples}.
}
\end{example}

{\begin{figure}[h!]
\centerline{\includegraphics{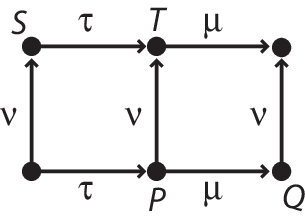}}
\caption{A medium on six states with $\bm=\mu\tilde\nu$ and $\bn=\tilde\nu\tau$.} \label{2squares} 
\end{figure}
}

\begin{corollary} \label{2.16}
The four following conditions are equivalent:

{\rm (i)} $ \ell(\bm) + \ell(\bn) \neq  \ell(\bm') + \ell(\bn')$ and $\mu\neq\tilde \tau$;

{\rm (ii)}  $\tau = \mu$;

{\rm (iii)} $\CCC(\bm) = \CCC(\bn)$;

{\rm (iv)}  $\ell(\bm) + \ell(\bn) + 2 =  \ell(\bm') + \ell(\bn')$.

Any of these conditions implies that $\tau\bm\tilde\mu\widetilde\bn$ is a $2$-gon for $S$. The converse does not hold.
\end{corollary}

\begin{proof}
We clearly have Case~6 of Theorem~\ref{6 cases}. Therefore, the four conditions are equivalent and $\tau\bm\tilde\mu\widetilde\bn$ is a $2$-gon for $S$. A counterexample is shown in Figure~\ref{2squares}.  It is easy to verify that the token system in Figure~\ref{2squares} is a medium. For $\bm=\mu\tilde\nu$ and $\bn=\tilde\nu\tau$ we have $S\tau\bm=S\bn\mu=Q$, so $\tau\bq\tilde\mu\widetilde\bw$ is a $2$-gon, but $\tau\not=\mu$. The medium in Figure~\ref{2squares} illustrates Case 3 of Theorem~\ref{6 cases}.
\end{proof}

\begin{theorem} \label{tokens<->semicubes}
The assignment $\tau\mapsto\WWW_\tau$ defines a one-to-one correspondence between $\TTT$ and the family of semicubes of the medium $(\SSS,\TTT)$.
\end{theorem}

\begin{proof}
Suppose that $\WWW_\tau=\WWW_\mu$ for some tokens $\tau,\mu\in\TTT$. There are states $S$, $T$, $P$, and $Q$ such that $S\tau=T$ and $P\mu=Q$, so $S,P\in\WWW_{\tilde\tau}=\WWW_{\tilde\mu}$ and $T,Q\in\WWW_\tau=\WWW_\mu$. There are three possible cases:

(i) $S=P$, $T\neq Q$. Let $\bm$ be a concise message transforming $Q$ into $T$. Since $T,Q\in\WWW_\tau$, neither $\tau$ not $\tilde\tau$ are in $\CCC(\bm)=\widehat{T}\setminus\widehat{Q}$ in contradiction to Lemma~\ref{small triangle}.

(ii) $S\neq P$, $T=Q$. Let $\bm$ be a concise message transforming $P$ into $S$. As in the previous case we have a contradiction with Lemma~\ref{small triangle}.

(iii) $S\neq P$, $T\neq Q$. As in case (i), neither $\tau$ not $\tilde\tau$ are in $\CCC(\bm)$ for a concise message $\bm$ transforming $Q$ into $T$. Thus we have Case 6 of Theorem~\ref{6 cases}. It follows that $\tau=\mu$.
\end{proof}

\section{Regular circuits} \label{S:regular}

Let $\bm=\tau_1\ldots\tau_{2n}$ be a closed message for a state $S$ of a medium $(\SSS,\TTT)$. (Note that, by Axiom [$\MMM$2], $\bm$ is vacuous and therefore must have an even length.) This closed message produces a cyclic sequence of states $(S_i)$ with $S_0=S_{2n}=S$. For any given $0\le i< 2n$, the message $\bm_i=\tau_{i+1}\ldots\tau_{2n}\ldots\tau_i$ is closed for $S_i$. (Note that $\bm_0=\bm$.) Thus we have $2n$ closed messages with the same cyclic sequence of states $(S_i)$.

{\begin{figure}[h!]
\centerline{\includegraphics{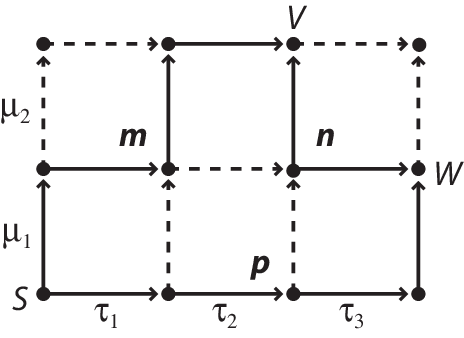}}
\caption{For Example~\ref{not orderly}.} \label{triangle3} 
\end{figure}
}

\begin{example} \label{not orderly}
{\rm The set of states of the token system in Figure~\ref{triangle3} are all vertices of the $4\times 3$-grid. 
Tokens indicated by parallel arrows are identical and only partly labeled. Reverse tokens are not shown. It can be verified that this token system is a medium. Let us define $\bm=\mu_1\tau_1\mu_2\tau_2$, $\bn=\tilde\mu_2\tau_3$, and $\bp=\tilde\mu_1\tilde\tau_3\tilde\tau_2\tilde\tau_1$. These messages are concise with $S\bm=V$, $V\bn=W$, and $W\bp=S$ (cf. Figure~\ref{triangle1}). The message $\bm\bn\bp$ is a closed message for $S$ of length $10$ indicated by bold arrows in Figure~\ref{triangle3}. As in the preceeding paragraph, one can produce nine more closed messages from $\bm\bn\bp$. It can be checked directly that none of these ten closed messages is a $2$-gon.
}
\end{example}

Example~\ref{not orderly} is an extreme instance of a closed message that does not define any $2$-gons. Another extreme case requires a definition. In what follows, $(S_i)$ is the sequence of states produced by a given message.

\begin{definition} \label{regular}
{\rm Let $\bm=\tau_1\ldots\tau_{2n}$ be a closed message for a state $S$. For $1\le i\le n$, the two tokens $\tau_i$ and $\tau_{i+n}$ are called {\em opposite}. The closed message $\bm$ is a {\em regular circuit} if the message $\tau_i\ldots\tau_{i+n-1}$ is concise for the state $S_{i-1}$ for all $1\le i\le n$.
}
\end{definition}

\begin{theorem} \label{opposite tokens}
Let $\bm=\tau_1\ldots\tau_{2n}$ be a $2$-gon for a state $S$. The following three conditions are equivalent:

{\rm(i)} $\bm$ is a regular circuit for $S$.

{\rm(ii)} For any $1\le i\le 2n-1$, the message $\tau_i\ldots\tau_{2n}\ldots\tau_{i-1}$ is a $2$-gon for the state $S_{i-1}$.

{\rm(iii)} The opposite tokens of $\bm$ are mutual reverses.
\end{theorem}

\begin{proof}
(i)~$\imp$~(ii). Since $\bm$ is a regular circuit, the message $\bp=\tau_i\ldots\tau_{i+n-1}$ is concise for the state $S_{i-1}$. The length of the message $\bq=\tau_{i+n}\ldots\tau_{2n}\ldots\tau_{i-1}$ is $n$, so we have $\dd(S_{i-1},S_{i+n-1})=\ell(\bp)=\ell(\bq)$. By Theorem~\ref{concise=shortest}, $\bq$ is a concise message for $S_{i+n-1}$. Thus $\bp\bq$ is a $2$-gon for $S_{i-1}$.

(ii)~$\imp$~(iii). Since the message $\tau_i\ldots\tau_{2n}\ldots\tau_{i-1}$ is a $2$-gon for the state $S_{i-1}$, the messages $\bp=\tau_i\ldots\tau_{i+n-1}$ and $\bq=\tau_{i+n}\ldots\tau_{2n}\ldots\tau_{i-1}$ are concise for the states $S_{i-1}$ and $S_{i+n-1}$, respectively. Thus, for $1\le i<n$ and 
$$
S=S_{i-1},\;T=S_i,\;Q=S_{i+n-1},\;P=S_{i+n},\;\tau=\tau_i,\;\mu=\tilde\tau_{i+n},
$$
we have Case~6 of Theorem~\ref{6 cases} (see Figure~\ref{rectangle2}). Hence, $\tau_i=\tilde\tau_{n+i}$.

(iii)~$\imp$~(i). Let $\bp=\tau_i\ldots\tau_{i+n-1}$ and $\bq=\tau_{i+n}\ldots\tau_{2n}\ldots\tau_{i-1}$. By Lemma~\ref{so-3-gon}, for any $1\le j\le n$, there is only one occurrence of the pair $\{\tau_j,\tilde\tau_j\}$ in the $2$-gon $\bm$. Since $\tilde\tau_j=\tau_{j+n}$ and $\ell(\bp)=\ell(\bq)=n$, there are no occurrences of $\{\tau_j,\tilde\tau_j\}$ in $\bp$ or $\bq$, so they are concise messages for $S_{i-1}$ and $S_{i+n-1}$, respectively. Thus $\bp\bq$ is a $2$-gon for $S_{i-1}$ and the result follows.
\end{proof}

As the following example illustrates, it is essential that the closed message $\bm$ in Theorem~\ref{opposite tokens} is a $2$-gon for some state.

\begin{example}
{\rm The message $\bm=\tau\tilde\nu\tilde\tau\nu\tau\tilde\nu\tilde\tau\nu\tau\tilde\nu\tilde\tau\nu$ is closed for the state $S$ of the medium shown in Figure~\ref{2squares}. The opposite tokens in $\bm$ are mutual reverses. Clearly, $\bm$ is not a regular circuit. Moreover, none of the $12$ circuits defined by $\bm$ for the produced states $S_i$ is a $2$-gon.
}
\end{example}

\section{Embeddings and isomorphisms} \label{S:embeddings}

The purpose of combinatorial media theory is to find and examine those properties of media that do not depend on a particular structure of individual states and tokens. For this purpose we introduce the concepts of embedding and isomorphism for token systems.

\begin{definition} \label{D:embedding}
{\rm Let $(\SSS,\TTT)$ and $(\SSS',\TTT')$ be two token systems. A pair $(\aa,\bb)$ of one--to--one functions $\aa:\SSS\rightarrow\SSS'$ and $\bb:\TTT\rightarrow\TTT'$ such that 
$$
S\tau=T \quad \eq \quad \aa\left(S\right)\bb\left(\tau\right)=\aa\left(T\right)
$$
for all $S,T\in\SSS$, $\tau\in\TTT$ is called an \emph{embedding} of the token system $(\SSS,\TTT)$ into the token system $(\SSS',\TTT')$.

Token systems $(\SSS,\TTT)$ and $(\SSS',\TTT')$ are \emph{isomorphic} if there is an embedding $(\aa,\bb)$ from $(\SSS,\TTT)$ into $(\SSS',\TTT')$ such that both $\aa$ and $\bb$ are bijections.
}
\end{definition}

Clearly, if one of two isomorphic token systems is a medium, then the other one is also a medium.

If a token system $(\SSS,\TTT)$ is a medium and $S\tau_1=S\tau_2\not=S$ for some state $S$, then, by Lemma~\ref{3 properties}(iii), $\tau_1=\tau_2$. In particular, if $(\aa,\bb)$ is an embedding of a medium into a medium, then $\bb(\tilde{\tau})=\widetilde{\bb(\tau)}$. Indeed, for a given $\tau$ there are two distinct states $S$ and $T$ such that $S\tilde\tau=T$. Then
\begin{align*}
\aa(S)\bb(\tilde\tau)=\aa(T)&\quad\eq\quad S\tilde\tau=T\quad\eq\quad T\tau=S\quad\eq\quad \\
	&\quad\eq\quad \aa(T)\bb(\tau)=\aa(S)\quad\eq\quad \aa(S)\widetilde{\bb(\tau)}=\aa(T),
\end{align*}
so $\bb(\tilde{\tau})=\widetilde{\bb(\tau)}$. We extend $\bb$ to the semigroup of messages  by defining 
$$
\bb(\tau_1\ldots\tau_k)=\bb(\tau_1)\ldots\bb(\tau_k).
$$
Clearly, the image $\bb(\bm)$ of a concise message $\bm$ for a state $S$ is a concise message for the state $\aa(S)$.

Let $(\SSS,\TTT)$ be a token system and $\QQQ$ be a subset of $\SSS$ consisting of more than two elements. The restriction of a token $\tau\in\TTT$ to $\QQQ$ is not necessarily a token on $\QQQ$. In order to construct a medium with the set of states $\QQQ$, we introduce the following concept.

\begin{definition} \label{D:restriction}
{\rm Let $(\SSS,\TTT)$ be a token system, $\QQQ$ be a nonempty subset of $\SSS$, and
$\tau\in\TTT$. We define a \emph{reduction} of $\tau$ to $\QQQ$ by 
$$
S\tau_{\QQQ} = \begin{cases}
	S\tau & \text{if $S\tau\in\QQQ$,} \\
	S & \text{if $S\tau\notin\QQQ$,}
\end{cases}
$$
for $S\in\QQQ$. A token system $(\QQQ,\TTT_{\QQQ})$ where
$\TTT_{\QQQ}=\{\tau_{\QQQ}\}_{\tau\in\TTT} \setminus\{\tau_0\}$ is the set of all distinct reductions of tokens in $\TTT$ to $\QQQ$ different from the identity function $\tau_0$ on $\QQQ$, is said to be the \emph{reduction} of $(\SSS,\TTT)$ to $\QQQ$.

We call $(\QQQ,\TTT_{\QQQ})$ a \emph{token subsystem} of $(\SSS,\TTT)$. If both $(\SSS,\TTT)$ and $(\QQQ,\TTT_{\QQQ})$ are media, we call
$(\QQQ,\TTT_{\QQQ})$ a \emph{submedium} of $(\SSS,\TTT)$.
}
\end{definition}

\begin{remark}
{\rm A reduction of a medium is not necessarily a submedium of a given medium. Consider, for instance, the medium shown in Figure~\ref{not submedium}.
The set of tokens of the reduction of this medium to $\QQQ=\{P,R\}$ is empty. Thus this reduction is not a medium.
}
\end{remark}

{\begin{figure}[h!]
\centerline{\includegraphics{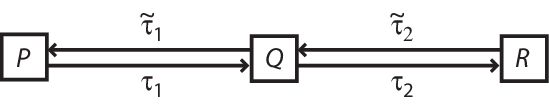}}
\caption{The reduction of this medium to $\{P,R\}$ is not a submedium.}\label{not submedium}
\end{figure}
}

The image $(\aa(\SSS),\bb(\TTT))$ of a token system $(\SSS,\TTT)$ under embedding $(\aa,\bb):(\SSS,\TTT)\rightarrow(\SSS',\TTT')$ is not, in general, the reduction of $(\SSS',\TTT')$ to $\aa(\SSS)$. Indeed, let $\SSS'=\SSS$, and let $\TTT$ be a proper nonempty subset of $\TTT'$. Then the image of $(\SSS,\TTT)$ under the identity embedding is not the reduction of $(\SSS,\TTT')$ to $\SSS$ (which is $(\SSS,\TTT')$ itself).

On the other hand, this is true in the case of media as the following theorem demonstrates.

\begin{theorem} \label{embedding theorem}
Let $(\aa,\bb):(\SSS,\TTT)\rightarrow(\SSS',\TTT')$ be an embedding of a medium $(\SSS,\TTT)$ into a medium $(\SSS',\TTT')$. Then the reduction $(\aa(\SSS),\TTT'_{\aa(\SSS)})$ is isomorphic to $(\SSS,\TTT)$.
\end{theorem}

\begin{proof}
For $\tau\in\TTT$, we define $\bb'(\tau)=\bb(\tau)_{\aa(\SSS)}$, the reduction of $\bb(\tau)$ to $\aa(\SSS)$. Let $S\tau=T$ for $S\not=T$ in $\SSS$. Then $\aa(S)\bb(\tau)=\aa(T)$ for $\aa(S)\not=\aa(T)$ in $\aa(\SSS)$. Hence, $\bb'$ maps $\TTT$ to $\TTT'_{\aa(\SSS)}$.

Let us show that $(\aa,\bb')$ is an isomorphism from $(\SSS,\TTT)$ onto $(\aa(\SSS),\TTT'_{\aa(\SSS)})$.

(i) $\bb'$ is onto. Suppose $\tau'_{\aa(\SSS)}\not=\tau_0$ for some $\tau'\in\TTT'$. Then there are $P\not= Q$ in $\SSS$ such that $\aa(P)\tau'_{\aa(\SSS)}=\aa(P)\tau'=\aa(Q)$. Let $Q=P\bm$ where $\bm$ is a concise message. We have
$$
\aa(Q)=\aa(P\bm)=\aa(P)\bb(\bm)=\aa(P)\tau',
$$
implying, by Theorem~\ref{Theorem1.16}, $\bb(\bm)=\tau'$, since $\bb(\bm)$ is a concise message. Hence, $\bm=\tau$ for some $\tau\in\TTT$. Thus $\bb(\tau)=\tau'$, which implies
$$
\bb'(\tau)=\bb(\tau)_{\aa(\SSS)}=\tau'_{\aa(\SSS)}.
$$
(ii) $\bb'$ is one--to--one. Suppose $\bb'(\tau_1)=\bb'(\tau_2)$. Since $\bb'(\tau_1)$ and $\bb'(\tau_2)$ are tokens on $\aa(\SSS)$ and $(\SSS',\TTT')$ is a medium, we have $\bb(\tau_1)=\bb(\tau_2)$. Hence, $\tau_1=\tau_2$.

(iii) Finally,
$$
S\tau=T \quad \eq \quad \aa\left(S\right)\bb'\left(\tau\right)=\aa\left(T\right),
$$
since
$$
S\tau=T \quad \eq \quad \aa\left(S\right)\bb\left(\tau\right)=\aa\left(T\right).
$$
This completes the proof.
\end{proof}

We conclude this section with an example of a submedium.

\begin{example} \label{submedium example}
{\rm Let $\FFF$ be a wg-family of finite subsets of a set $X$. The representing medium $(\FFF,\GGG_\FFF)$ of $\FFF$ is clearly the reduction of the complete medium $(\BBB(X),\GGG_{\BBB(X)})$ to $\FFF$. Thus, $(\FFF,\GGG_\FFF)$ is a submedium of $(\BBB(X),\GGG_{\BBB(X)})$ for any wg-family $\FFF$.
}
\end{example}

\section{Media and partial cubes} \label{S:p-cubes}

Let $G$ be the graph of a medium $(\SSS,\TTT)$. As we observed before, $G$ is a connected bipartite graph, but not any connected bipartite graph is a graph of a medium (see Example~\ref{bipartite not medium}). By Theorems~\ref{convex semicubes},~\ref{metric semicubes}, and~\ref{DWT}, the graph $G$ is a partial cube:

\begin{theorem} \label{media->cube}
The graph $G$ of a medium $(\SSS,\TTT)$ is a partial cube.
\end{theorem}

We give two more proofs of this important result. The first proof utilizes the concept of Winkler's relation $\Theta$ (see~\ref{S:bipartite}).

\begin{proof}
The four vertices $S$, $T$, $P$, and $Q$ of two edges $\{S,T\}$ and $\{P,Q\}$ that stand in Winkler's relation $\Theta$ form a `rectangle' described by Cases 5 and 6 of Theorem~\ref{tetrahedron theorem} (see Figure~\ref{cases5-6} and the text following that figure). Since edges of the graph $G$ correspond to pairs of mutually reverse tokens, it follows from Theorem~\ref{6 cases} that two edges of $G$ stand in the relation $\Theta$ if and only if they represent the same pair of mutually reversed tokens. Thus the relation $\Theta$ is an equivalence relation on the set of edges of $G$. By Theorem~\ref{DWT}, $G$ is a partial cube.
\end{proof}

The second proof is based on the result of Theorem~\ref{labels}.

\begin{proof}
The edges of $G$ are labeled by elements of the set $J=\{\{\tau,\tilde\tau\}\}_{\tau\in\TTT}$. Since the shortest paths of $G$ correspond to the concise messages of $(\SSS,\TTT)$, condition (i) of Theorem~\ref{labels} is satisfied. A closed walk $W$ in $G$ defines a closed message $\bm$ for a vertex of $W$. By Axiom [$\MMM$2], the message $\bm$ is vacuous. Thus every label appears an even number of times in the walk $W$. The result follows from Theorem~\ref{labels}.
\end{proof}

Let $(\SSS,\TTT)$ be a medium and $G$ be its graph. By Theorem~\ref{media->cube}, $G$ is a partial cube, so there is an isometric embedding $\aa$ of $G$ into the cube $\HHH(X)$ for some set $X$. The set $\aa(\SSS)$ is a wg-family $\FFF$ of finite subsets of $X$. Let $(\FFF,\GGG_\FFF)$ be the representing medium of this wg-family. These objects are schematically shown in the diagram below, where $\langle\FFF\rangle$ is an isometric subgraph of $\HHH(X)$ induced by the family $\FFF$.
\begeq \label{medium<->rep.medium}
(\SSS,\TTT)\quad\xrightarrow[\text{of $(\SSS,\TTT)$}]{\text{graph}}\quad G\quad\xrightarrow{\;\aa\;}\quad\langle\FFF\rangle\quad\xrightarrow[\text{medium}]{\text{representing}}\quad(\FFF,\GGG_\FFF)
\edeq

\begin{theorem} \label{isomorphism theorem}
The media $(\SSS,\TTT)$ and $(\FFF,\GGG_\FFF)$ are isomorphic.
\end{theorem}

\begin{proof}
Clearly, $\aa$ is a bijection from $\SSS$ onto $\FFF$. By Theorem~\ref{metric semicubes} and Definition~\ref{semicube graph}, the semicubes of $(\SSS,\TTT)$ and $G$ (resp. $(\FFF,\GGG_\FFF)$ and $\langle\FFF\rangle$) are defined in terms of the metric $\dd$ on $\SSS$ (resp. the metric $d$ on $\FFF$). Since $\aa$ is an isometric embedding, it defines a one-to-one correspondence between semicubes of $G$ and $\langle\FFF\rangle$. By Theorem~\ref{tokens<->semicubes}, $\tau\mapsto\WWW_\tau$ (resp. $\gg\mapsto\WWW_\gg$) is a bijection from $\TTT$ (resp. $\GGG_\FFF$) onto the family of semicubes of the medium $(\SSS,\TTT)$ (resp. $(\FFF,\GGG_\FFF)$). The above bijections define a bijection $\bb:\TTT\rightarrow\GGG_\FFF$:
$$
\bb:\quad\TTT\quad\rightarrow\quad\{\WWW_\tau\}_{\tau\in\TTT}\quad\xrightarrow{\aa}\quad\{\WWW_\gamma\}_{\gg\in\GGG_\FFF}\quad\rightarrow\quad\GGG_\FFF
$$
Since a pair of opposite semicubes form a partition of the set of states of a medium, we have $\bb(\tilde\tau)=\widetilde{\bb(\tau)}$. 
Suppose that $S\tau=T$ for $S\neq T$. Then $S\in\WWW_{\tilde\tau}$ and $T\in\WWW_\tau$. Therefore, 
$$
\aa(S)\in\WWW_{\bb(\tilde\tau)}=\WWW_{\widetilde{\bb(\tau)}}\quad\text{and}\quad \aa(T)\in\WWW_{\bb(\tau)},
$$
so $\aa(S)\bb(\tau)=\aa(T)$. Hence, 
$$
S\tau=T\quad\imp\quad\aa(S)\bb(\tau)=\aa(T),
$$
for $S\neq T$. A similar argument shows that the converse is also true. It follows that $(\aa,\bb)$ is an isomorphism from $(\SSS,\TTT)$ onto $(\FFF,\GGG_\FFF)$.
\end{proof}

\begin{remark}
{\rm The representing medium $(\FFF,\GGG_\FFF)$ in~(\ref{medium<->rep.medium}) is not defined uniquely because there are many possible embedding $\aa$ of the partial cube $G$ into a cube. Each of these embeddings defines a particular wg-family $\FFF$ and these families are quite different. On the other hand, as it follows from Theorem~\ref{isomorphism theorem}, all representing media defined by~(\ref{medium<->rep.medium}) are isomorphic. The correspondence
$$
\text{medium $(\SSS,\TTT)$}\quad\mapsto\quad\text{graph of $(\SSS,\TTT)$}
$$
defines a bijection from the set of media on the set of states $\SSS$ onto the set of partial cubes with the vertex set $\SSS$.
}
\end{remark}

\begin{remark}
{\rm If the graph of a token system is a partial cube, it does not mean necessarily that the token system itself is a medium (see Figure~\ref{not medium}).
}
\end{remark}

{\begin{figure}[h!]
\centerline{\includegraphics{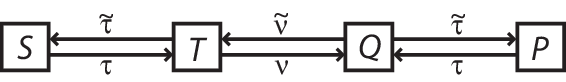}}
\caption{This token system is not a medium since there is no concise message producing the state $P$ from the state $S$. Its graph is a partial cube.}\label{not medium}
\end{figure}
}

As Winkler's relation $\Theta$ combines Cases 5 and 6 of Theorem~\ref{tetrahedron theorem} (see the second proof of Theorem~\ref{media->cube}), the relation $\bmfL$ is defined on the set of arcs (ordered pairs of adjacent vertices) of a graph by using equations~(\ref{case6}) from Case~6 of Theorem~\ref{tetrahedron theorem}:
$$
(S,T)\bmfL(P,Q)\quad\eq\quad \quad\dd(S,P)=\dd(T,Q)=\dd(T,P)-1=\dd(S,Q)-1.
$$
By Theorem~\ref{6 cases} (Case~6), the relation $\bmfL$ is transitive on the set of arcs of the graph of a medium. Clearly, $\bmfL$ is reflexive and symmetric, so it is an equivalence relation in the case of the graph of a medium. Connected bipartite graphs for which $\bmfL$ is an equivalence relation are called {\em mediatic} in~\cite{dE07,jF07}, so we have the following result.

\begin{theorem} \label{media->mediatic}
The graph of a medium is mediatic.
\end{theorem}

It can be shown that the class of partial cubes coincides with the class of mediatic graphs.

\begin{remark}
{\rm Note that, unlike in Case 6, the equations in~(\ref{case5}) (Case 5) do not define a transitive relation on the set of arcs of $G$. This is illustrated by the drawing in Figure~\ref{nontransitive}.
}
\end{remark}

{\begin{figure}[h!]
\centerline{\includegraphics{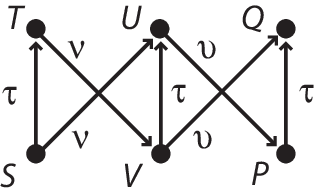}}
\caption{Intransitivity of the relation defined by equations in~(\ref{case5}).} \label{nontransitive}
\end{figure}
}

\section*{Acknowledgments}

The author is greatly indebted to Jean-Claude Falmagne and David Eppstein for many generous discussions on media theory.

\end{document}